\documentclass{article}
\usepackage{latexsym}
\usepackage{amsthm}  
\usepackage{amssymb}
\usepackage{amsmath}

\newcommand{\mbb}[1]{\mathbf{#1}}

\newcommand{\msc}[1]{\mathcal{#1}}
\newcommand{\mrm}[1]{\mathrm{#1}}
\newcommand{\Pro}{\mbb{P}}
\newcommand{\sh}[1]{\msc{#1}}

\newcommand{\pd}[2]{{\frac{\partial{#1}}{\partial{#2}}}}
\newcommand{\ses}[5]{0\rightarrow{#1}\stackrel{{#2}}{\rightarrow}{#3}
\stackrel{{#4}}{\rightarrow}{#5}\rightarrow 0}

\theoremstyle{definition}
\newtheorem{defn}{Definition}
\newtheorem{rmk}[defn]{Remark}

\newtheorem{exa}[defn]{Example}
\theoremstyle{plain}
\newtheorem{thm}[defn]{Theorem}

\newtheorem{lem}[defn]{Lemma}
\newtheorem{prop}[defn]{Proposition}

\date{}		

\begin{document}

\title{Variation of hyperplane sections}
\author{Michael A. van Opstall, R\u{a}zvan Veliche}
\maketitle

Throughout this article we assume $n>2$ and $\epsilon^2=0$. All varieties
are defined over an algebraically closed field $k$.

In \cite{hmp:hld}, Harris, Mazur, and Pandharipande ask: given a
smooth, degree $d$ hypersurface $X$ in $\Pro^n$, does the family of
all smooth $l$-plane sections of $X$ vary maximally in moduli? They
answer this question when $d$ and $l$ are quite small compared to $n$, and
when the characteristic of the base field is zero. The
purpose of this article is to establish that the variation of smooth
hyperplane sections of $X$ is maximal when $X$ is a general smooth
hypersurface of degree $d$, over an algebraically closed field of arbitrary
characteristic.

Beauville proves that the variation in moduli is not zero (as long as
$d>2$, of course), except for Fermat hypersurfaces of degree $d$ such that
$(d-1)$ is a power of the characteristic \cite{b:vm}. The proof is
elegant, and we follow the same line here, at least initially. We also show
that the variation of hyperplane sections of the Fermat is maximal except
in those cases excluded by Beauville.

The main theorem of this article is

\begin{thm}\label{mainthm}
If $X$ is a general hypersurface in $\Pro^n$, then the
hyperplane sections of $X$ vary maximally in moduli.
\end{thm}

First of all, note that the degree one and two cases are trivial, since
there are no moduli, so we assume from now on that $d>2$. Geometric
invariant theory then provides a
moduli space $M(d,n)$ of smooth degree $d$ hypersurfaces in $\Pro^n$, at least
when $\mrm{char~}k\not | ~d$ (since the methods employed in this paper are
infinitesimal, we need not rely on the existence of a global moduli space, but
we use it for motivation).
This is the quotient of some open set of the projective space
parameterizing all degree $d$ hypersurfaces by the action of
$\mrm{Psl}(n+1)$, so the dimension $m(d,n)$ of the moduli space is
\[
m(d,n)=\binom{n+d}{d}-(n+1)^2.
\]
The space of hyperplane sections is $n$ dimensional, so $n>m(d,n-1)$ if
and only if $n=3$ and $d=3$. That is, for a cubic surface, there is
a three parameter family of hyperplane sections, but only one modulus
of cubic curves. In this case we have

\begin{thm}
Suppose $\mrm{char~}k\neq 2$, and let $X$ be a smooth cubic surface. Then the
smooth hyperplane sections have maximal variation in moduli. If $\mrm{char~}k$
is two, then if $X$ is a general cubic surface the smooth hyperplane
sections have maximal variation in moduli.
\end{thm}

\begin{proof}
This follows from the main result of Beauville quoted above.
\end{proof}

Therefore, in what follows, we may assume that $n>2$, $d>2$, and
$(n,d)\neq(3,3)$. Then we are asking that the map from the complement
of the discriminant locus in the dual projective space parameterizing
hyperplane sections of $X$ to the moduli space of hypersurfaces is
generically finite. Note that this is in a sense the opposite of the
considerations of Harris, Mazur, and Pandharipande, in that their lower
bound for $n$ in terms of $d$ and $k$ ensures that maximal variation is
checked by checking {\em surjectivity} of the differential of the variation,
whereas we will check {\em injectivity}.

\section{Infinitesimal study}

In this section, we reduce the problem to linear algebra by considering
it only to first order.

\subsection{Preliminaries}

As above, let $X$ be a smooth degree $d$ hypersurface in $\Pro^n$. Let
us briefly recall some basics about deformations of $X$. Since $X$ is
smooth, its first order deformations up to isomorphism are classified by
$H^1(X,T_X)$. The
deformations of $X$ as a subvariety of $\Pro^n$ are classified by
$H^0(X,N_{X/\Pro^n})$. In this case, $N_{X/\Pro^n}\cong \sh{O}_X(d)$,
and choosing a degree $d$ polynomial $g$ in this space, the corresponding
first order deformation is given by $f+\epsilon g=0$. The short exact
sequence
\[
\ses{T_X}{}{T_{\Pro^n}|_X}{}{N_{X/\Pro^n}}
\]
induces a morphism $H^0(N_{X/\Pro^n})\rightarrow H^1(T_X)$ taking an
embedded deformation to its isomorphism class as an abstract
deformation. Recall that the Jacobian ring $R(X)$ of a hypersurface $X$
with defining equation $f$ is the ring
\[
k[x_0,\ldots,x_n]/(f,\partial f/\partial x_0, \ldots, \partial f/\partial x_n).
\]
Note that if $\mrm{char~}k \not |~d$, then $f$ is automatically
in the ideal generated by the partials by the Euler formula. We will
use the following result of Beauville ({\em loc. cit.}), whose proof is
elementary:

\begin{lem}
The morphism $H^0(N_{X/\Pro^n})\rightarrow H^1(T_X)$ factors through
$R(X)_{d}$, and the induced morphism $R(X)_d\rightarrow H^1(T_X)$ is
injective.
\end{lem}

\begin{rmk}
This lemma is false in the case $n=2$ and $d=3$, which is excluded
from our consideration.
\end{rmk}

\subsection{Criterion for maximal variation}

Now let $H$ be a
hyperplane such that $X\cap H$ is smooth. Suppose $f$ is the defining
equation of $X$ and that $H$ is given by $x_0=0$. Then to first order,
a deformation of $H$ is given by $x_0=\epsilon l(x_1,\ldots,x_n)$ where
$l$ is a linear form. The equation for the corresponding
first order deformation of $X\cap H$ is given by
\[
f(\epsilon l(x_1,\ldots,x_n),x_1,\ldots,x_n)=0.
\]
We may expand this in $\epsilon$ to obtain
\[
f(0,x_1,\ldots,x_n)+\epsilon\pd{f}{x_0}(0,x_1,\ldots,x_n)l(x_1,\ldots,x_n)
=0
\]

From the results of the last section, we conclude:

\begin{prop}\label{crit}
Notation as above. Suppose $\pd{f}{x_0}(0,x_1,\ldots,x_n)\neq 0$.
An embedded deformation of $X\cap H$ corresponding to a linear form
$l$ in $x_1,\ldots, x_n$ is trivial to first order (as an abstract
deformation) if and only if
\[
\pd{f}{x_0}(0,x_1,\ldots,x_n)l(x_1,\ldots,x_n)
\]
is zero in the
Jacobian ring $R(X\cap H)$.
\end{prop}

\begin{rmk}
Note that this first order criterion must be applied with care.
For example, for the hyperplane section
$x_0=0$ of the Fermat hypersurface, the left hand side will be zero regardless
of the choice of $l$. On the other hand, in general, the
hyperplane sections of the Fermat have some variation (and in characteristic
zero, in fact, maximal variation; see the examples below). We must show that
$l$ vanishes, so choosing a hyperplane section as above with
$\pd{f}{x_0}(0,x_1,\ldots,x_n)=0$ gives us no information. When this
derivative vanishes, we must perturb our hyperplane section a little and then
check the criterion.
\end{rmk}

\section{Openness of maximal variation}

\begin{thm}\label{open}
The set of smooth hypersurfaces whose variation of hyperplane sections
is maximal is Zariski open in the space of all smooth hypersurfaces.
\end{thm}

\begin{proof}
Denote by $D(d,n)$ the open subset of projective space parameterizing smooth
hypersurfaces of degree $d$ in $\Pro^n$.
Let $U$ be the set in $D(d,n)$ whose hyperplane sections
vary maximally in moduli. Let $X$ be such a hypersurface. Then there
exists a hyperplane section $H$ such that the map
\[
H^0(N_{X\cap H/H})\rightarrow H^1(T_{X\cap H})
\]
is injective. But $H^0(N_{X\cap H/H})$ can be identified with the tangent
space to $D(d,n-1)$ at the point corresponding to $X\cap H$. Let $\sh{X}$ be
the universal hypersurface over $D(d,n)$ and $\sh{H}$ the constant
hyperplane $D(d,n)\times H$. Let $\sh{Y}$ be the intersection
$\sh{X}\cap\sh{H}$ and $\pi$ the projection of $\sh{Y}$ onto
$D(d,n-1)$. Note that since $X$ will have singular hyperplane sections,
$\pi$ is not everywhere defined, but is defined in a Zariski open
neighborhood of $X\cap H$. Then the map of cohomology spaces above
is just the restriction of the map of vector bundles
\[
T_{D(d,n-1)}\rightarrow R^1\pi_*T_{\sh{Y}/D(d,n-1)}
\]
to the point in $D(d,n-1)$ corresponding to $X$. Since this map of
vector bundles is injective at this point, it is injective in an open
neighborhood. Hence $U$ is open.
\end{proof}

Theorem \ref{mainthm} follows from this theorem if the set $U$ is nonempty.
This is shown in the next section.

\section{Examples}

\begin{exa}
The hyperplane sections of the smooth hypersurface defined by
\[
\sum_{i=0}^n x_i^d+\sum_{j=0}^{n-1}x_j^{d-1}x_{j+1}+x_n^{d-1}x_0=0
\]
vary maximally in moduli.
\end{exa}

\begin{proof}
We will use the criterion from Proposition \ref{crit}. Small deformations
of the hyperplane section $x_0=0$ vary maximally in moduli if when we write
\begin{equation}\label{eqn}
l_0x_n^{d-1}=\sum_{i=1}^{n-1}l_i(dx_i^{d-1}+(d-1)x_i^{d-2}x_{i+1}
+x_{i-1}^{d-1})+l_n(dx_n^{d-1}+x_{n-1}^{d-1}),
\end{equation}
where the $l_i$ are linear forms in the variables $x_1,\ldots,x_n$, we
can conclude that $l_0=0$.

First assume $\mrm{char~}k ~|~(d-1)$. Then by considering terms in which
$x_n$ occurs to the $(d-1)$st or $d$th power, we conclude $l_0=dl_n$.
Repeating for terms where $x_{n-1}$ occurs to the $(d-1)$st or $d$th power
we obtain $l_n=-dl_{n-1}$. Continuing likewise, we see that all the $l_i$
are multiples of $l_1$.
We can
cancel $l_1$ from the rewritten form of (\ref{eqn}). In this way, we obtain
a linear relation among the partial derivatives of $f$ with $x_0$ set equal
to zero, which contradicts the smoothness of $X$.

The second case is when $\mrm{char~}k ~|~d$. Considering terms with high
powers of $x_n$ as above, if $d>3$ it follows immediately that $l_0=0$.
If $d=3$, we see that necessarily
\[
l_{n-1}=ax_n+bx_{n-1}+cx_{n-2}, l_0=2ax_{n-1},
\]
but this introduces a term on the right hand side $ax_nx_{n-2}^2$ which
cannot be cancelled by any other term, so $a=0$.

So we may assume that $\mrm{char~}k$ divides neither $(d-1)$ nor $d$. Again,
for simplicity, assume first that $d>3$. Then consindering terms with high
powers of $x_n$, we see that $l_0=dl_n$. Passing on to terms with $x_{n-1}$
occuring to the power of at least $d-2$, cancelling $x_{n-1}^{d-2}$ we obtain
\[
0=dl_{n-1}x_{n-1}+(d-1)l_{n-1}x_n+l_nx_{n-1}
\]
from which it follows that
\[
l_{n-1}=ax_{n-1}
\]
and that $l_n=0$ (and therefore $l_0=0$) if $l_{n-1}=0$. Considering terms
with $x_{n-2}$ occuring to a high power, we get
\[
0=dl_{n-2}x_{n-2}+(d-1)l_{n-2}x_{n-1}+ax_{n-1}x_{n-2}
\]
and conclude that $l_{n-2}=\alpha x_{n-2}+\beta x_{n-1}$ for
some choice of constants. Plugging this back in shows that in fact
$l_{n-2}$ has to be zero, which also implies that $l_{n-1}$ is zero, so
we are done.

The case where $d=3$ is similar and left to the reader.
\end{proof}

The second example is superfluous in proving the main result, but shows that
in some sense the ``most probable counterexample'' to the conjecture that
hyperplane sections vary maximally (at least in ``good'' characteristics) is
in fact not a counterexample.

\begin{exa}\label{fermat}
Assume $\mrm{char~}k\neq d$ and that $d-1$ is not a power of
$\mrm{char~}k$. Then the hyperplane sections of
the Fermat hypersurface defined by $\sum x_i^d=0$ vary maximally in moduli.
\end{exa}

\begin{proof}
Since $\mrm{char~}k\neq d$, the Fermat is smooth.
As noted above, the hypersurface section defined by $x_0=0$ is not a good
choice for applying our criterion. Let
$a=\sum_{i=1}^n a_ix_i$ be a general linear form. Then we apply our
criterion to the hypersurface defined by
\[
(x_0+a)^d+\sum_{i=1}^n x_i^d=0,
\]
which is equivalent to considering variation near the hypersurface
section $x_0+a=0$ of the Fermat (which is smooth by Bertini, since $a$ is
general). We assume there is a relation of the form:
\[
l_0\pd{f}{x_0}\left|_{x_0=0}\right.=\sum_{i=1}^n l_i\pd{f}{x_i}\left|_{x_0=0}
\right.
\]
where $l_i$ are linear forms in $x_1,\ldots, x_n$ as above. That is,
\[
dl_0a^{d-1}=\sum dl_i(a_ia^{d-1}+x_i^{d-1}).
\]
Set $m=l_0-\sum_{i=1}^n a_il_i$, so that
\[
ma^{d-1}=\sum l_i x_i^{d-1}.
\]
Since $a$ is a general linear form, the polynomial on the left hand side
has monomials with three or more distinct variables (as long as $d-1$ is
not a power of the characteristic, so ``freshman exponentiation'' does not
hold), but the right hand
side does not, so $m=0$, from which it follows that all the $l_i$ are zero.
\end{proof}

\section{Conclusion}

Consider the hypersurface defined by
\[
0=x_0^3+x_1^3+x_0x_1^2+x_1x_2^2+x_3^3+x_2x_4^2.
\]
It is smooth, and the hyperplane section $x_0=0$ is also smooth. Furthermore,
since $\pd{f}{x_0}\left|_{x_0=0}\right.\neq 0$, the criterion is not
vacuous. However,
if $l_0=ax_1+bx_4$, one can solve for $l_1,\ldots,l_4$ in the criterion
above. So the variation of hyperplane sections is not maximal near this
hyperplane. However, one can check that variation is still maximal near some
other hyperplane. Here, the two tangent vectors which are killed must be
tangent to some two-dimensional subvariety which is blown down by the map to
the moduli space. So one must check the criterion for all possible
linear perturbations of $x_0$. This is computationally quite complex, since
the computations must be done symbolically, and
at present, the computation for sections of a cubic threefold seems too
intense for Maple (at least for the authors' patience), even when the form of
the equation can be simplified using
coordinate changes in certain characteristics. For example, in characteristic
zero, the equation of a general cubic threefold can be written in the form
\[
x_0^3+x_0\left(\sum_{i=1}^4a_ix_i^2\right)+g(x_1,x_2,x_3,x_4)
\]
where $g$ is a cubic form in four variables and the $a_i$ are constant.

\end{document}